\documentclass[12pt]{article}
\usepackage{textcomp}
\usepackage{url}
\usepackage{amsthm}
\usepackage{amsmath}
\usepackage{fullpage}
\usepackage{enumerate}
\usepackage{amssymb}

\usepackage{booktabs}
\usepackage{float}

\newtheorem{thm}{Theorem}

\newtheorem{prob}{Problem}

\title{Uchiyama's conjecture on sums of squares}
\author{Tim Trudgian\footnote{Supported by Australian Research Council Future Fellowship FT160100094.}\\
 {\small School of Physical, Environmental and Mathematical Sciences}\\ {\small  The University of New South Wales Canberra, Australia }\\
 {\small  t.trudgian@adfa.edu.au }}

\begin{document}
\maketitle
\noindent\textit{To the Australian Test Cricket selectors, to whom I am indebted a slice of humble pie, after they ignored the strident remonstrances of TS Trudgian and included MR Marsh in the team.}
\begin{abstract}
\noindent
Uchiyama showed that every interval $(n, n + c n^{1/4})$ contains an integer that is the sum of two squares, where $c= 2^{3/2}$. He also conjectured a minimal value of $c$ such that the above statement still holds. We investigate this claim.
  \end{abstract}

\section{Introduction}
We are in the market for short intervals $(x, x+h(x)]$ where $h(x) = o(x)$ such that each interval (or at least each interval with $x$ sufficiently large) contains a number that is the sum of two squares. The analogous question for primes in short intervals is well-known. The best result is $h(x)= x^{0.525}$ by Baker, Harman and Pintz \cite{BHP}, which is true for sufficiently large\footnote{Indeed, this result is effective in the sense that one \textit{could} calculate the $x_{0}$ --- presumably astronomically large --- such that the result holds for all $x\geq x_{0}$. Setting this problem to a graduate student could well violate the 8th Amendment to the US Constitution.} $x$. 

Anything that can be done for primes can presumably be done for the sum of squares. After all, every prime $p\equiv 1 \pmod{4}$ is the sum of two squares. 
Indeed, we have that
\begin{equation}\label{jones}
\pi(x) \sim \frac{x}{\log x}, \quad R(x) \sim C\frac{x}{\sqrt{\log x}},
\end{equation}
where $\pi(x)$ denotes the number of primes not exceeding $x$, and $R(x)$ the number of $n\leq x$ that are the sum of two squares. Here $C$ is a constant $\approx 0.74$ the exact form of which is known but does not concern us here. Therefore, by (\ref{jones}) there are many more sums of squares than there are primes.

Perhaps, then, it is not surprising  that one can do better in the sum-of-squares case. Indeed, one can show easily that $h= cx^{1/4}$ is admissible, for some constant $c$ --- see, e.g., Heath-Brown \cite[p.\ 2]{HB}. Bambah and Chowla \cite{Bambah} proved that one may take $c= 2^{3/2} + \epsilon$ where $\epsilon\rightarrow 0$ as $x\rightarrow \infty$. Uchiyama \cite{U} neatened this up and showed (in less than two pages, using nothing more than introductory calculus) that one could take $c = 2^{3/2}$ for all $x$. Uchiyama's proof is generalised to other quadratic forms  in \cite{Kaplan}. 
Remarkably, the elementary methods used in these papers have not been improved upon. Therefore, one may pose the following.
\begin{prob}\label{P1}
Find a value of $c< 2^{3/2}$ such that for all $n\geq 1$ (or even all $n$ sufficiently large) there are integers $x,y$ such that
\begin{equation}\label{indiana}
n< x^{2} + y^{2} < n + c n^{1/4}.
\end{equation}
\end{prob}
At the 61st Annual Meeting of the Australian Mathematical Society, held recently at Macquarie University, I offered a bottle of wine (or cricket merchandise vouchers for non-drinkers) to the tune of \$50 AUD for a resolution of Problem \ref{P1}. Freed from the pressure of prizes, one may also wish to consider the following.
\begin{prob}\label{P2}
Assume the Riemann hypothesis. Find a function $h(x) = o(x^{1/4})$  such that for all $x$ sufficiently large there is  a sum of squares in  $(x, x+ h(x)]$.
\end{prob}
Problem \ref{P2} is presumably much harder than Problem \ref{P1}  and hence deserves a larger prize for its resolution. Unfortunately the margins in my chequebook are too small to write down a full amount.

Again, on analogy with primes in short intervals, one may examine Cram\'{e}r's conjecture, viz.\ that one may take $h= c \log^{2} x$ for some constant $c$. It is known that $h(x)$ cannot be smaller than a constant times 
$$\frac{(\log x) (\log\log x) (\log\log\log\log x)}{\log\log\log x},$$
see --- \cite{Ford} and the references therein for an overview of historical developments.

To this end, and keeping in mind that anything primes can do sums-of-squares ought to do better, Erd\H{o}s conjectured (see \cite[\S 18]{Erdos}) that one could take $h= x^{\epsilon}$. He also showed \cite{Erdos50} that $h(x)$ could not be made smaller than a constant times $\log x/(\log\log x)^{1/2}$. 
\begin{prob}\label{P3}
Assume the Riemann hypothesis (and any other high-powered conjecture you like). Give a conjecture on the smallest $h(x)$ such that for all $x$ sufficiently large there is a sum of squares in $(x, x+ h(x)]$.
\end{prob}
Combining Erd\H{o}s' result with Cram\'{e}r's conjecture, it seems reasonable that one should look for a function $h(x)$ in Problem \ref{P3} such that
$$ \frac{\log x}{(\log\log x)^{1/2}} \ll h(x) \ll \log^{2} x.$$

\section{Uchiyama's conjecture}
We now return to Problem \ref{P1}. Uchiyama shows that were one looking to improve the constant $c$ in (\ref{indiana}) for all $n\geq 1$, one must have $c>2^{-1/2} 5^{3/4} = 2.364\ldots$. This comes about from inserting $n=20$ into (\ref{indiana}) and noting that none of $21, 22, 23, 24$ is a sum of two squares. Uchiyama conjectures that this is the worst value of $c$, that is,  (\ref{indiana}) should be true for all $c>2^{-1/2} 5^{3/4}$. It is remarked on \cite[p.\ 126]{U} that this has been verified for all $n\leq 1000$.

One can easily extend this and, indeed, note that at $n=1493$ we run into a problem. Since none of $1494, 1495, \ldots, 1507$ is a sum of two squares, we find that if (\ref{indiana}) is true then we must have $c> 15/(1493)^{1/4} = 2.413\ldots$. We find that this repaired conjecture of Uchiyama's is valid for all $n\leq 10^{8}$. To summarise, we have
\begin{thm}\label{TT}
There is a constant $c\in(15/(1493)^{1/4}, 2^{3/2}] = (2.413\ldots, 2.828\ldots]$ such that for all $n\geq 1$ there are integers $x,y$ such that
\begin{equation}\label{indiana2}
n< x^{2} + y^{2} < n + c n^{1/4}.
\end{equation}
Moreover, for any constant $c>15/(1493)^{1/4}$ the inequality (\ref{indiana2}) is true for all $n\leq 10^{8}$.
\end{thm}
It should not be difficult to extend the range $10^{8}$ in Theorem \ref{TT}. Indeed, the code used was slovenly in the extreme. We merely enumerated all intervals in (\ref{indiana2}): once we found a value inside an interval that was the sum of two squares, we ticked off this $n$, and went to the next one. This took  6 hours using \textit{Mathematica} on a standard desktop machine. This overlooks entirely the fact that checking one interval will suffice for many others. For example, with $n=400$ we have the interval $(401, 410.79\ldots)$. Since $410$ is a sum of two squares this would serve to verify (\ref{indiana}) for $n=400, 401, \ldots, 409$ --- ten for the price of one!

\subsubsection*{Acknowledgements}
I wish to thank Xuan Duong and Paul Smith for organising a cracking AustMS meeting at Macquarie University, where I presented this material. Thanks also to Dzmitry Badziahin, Kieran Clancy, Adrian Dudek, Mike Mossinghoff and Gerry Myerson for rollicking discussions on the topic.

\end{document}